# Decomposition of a Nonlinear Multivariate Function using the Heaviside Step Function

## Eisuke Chikayama[1,2,3,*]

*1 Department of Information Systems, Niigata University of International and Information Studies, 3-1-1 Mizukino, Nishi-ku, Niigata-shi, Niigata 950-2292, Japan*

*2 Environmental Metabolic Analysis Research Team, RIKEN, 1-7-22 Suehiro-cho, Tsurumi-ku, Yokohama-shi, Kanagawa 230-0045, Japan*

*3 Image Processing Research Team, RIKEN, 2-1 Hirosawa, Wako-shi, Saitama 351-0198, Japan*

* Author to whom correspondence should be addressed. Tel.: +81-25-239-3706. Fax: +81-25-239-3690. E-mail: chikaya@nuis.ac.jp.

**Abstract**

Whereas the Dirac delta function introduced by P. A. M. Dirac in 1930 in his famous quantum mechanics text has been well studied, a not famous formula related to the delta function using the Heaviside step function in a single-variable form, also given in Dirac's text, has been poorly studied. We demonstrate the decomposition of a nonlinear multivariate function into a sum of integrals in which each integrand is composed of a derivative of the function and a direct product of Heaviside step functions. It is an extension of Dirac's single-variable form to that for multiple variables. Moreover, it remains mathematically equivalent to the definition of the Dirac delta function with multiple variables, and offers a mathematically unified expression.

## 1. Introduction

P. A. M. Dirac introduced in 1930 a function, now called the Dirac delta function, to develop his theory of quantum mechanics [1]. The delta function, a valuable function in the field of mathematical physics, takes value infinity at $x = 0$ and zero at $x \neq 0$; its integral is unity. Its fundamental property derivable from its definition is that any nonlinear multivariate real function can be expressed with delta functions $\delta$ and integrals as follows,

$$R(X_1, X_2, \cdots, X_N) = \int_{-\infty}^{\infty} \cdots \int_{-\infty}^{\infty} R(\mu_1, \cdots, \mu_N)\delta(\mu_1 - X_1) \cdots \delta(\mu_N - X_N) d\mu_1 \cdots d\mu_N \,. \tag{1.1}$$

The importance of this property is analogous to the Fourier transform [2] for its ability to yield an alternative representation of any nonlinear multivariate function. Moreover, the delta function can be applied as an alternative to the Kronecker delta from the view point of continuity, and applied to numbers of formulae in the Fourier and Laplace transforms, and differential equations [3]. A more rigorous mathematical theory for the delta function has also been developed and expanded as the theory of distributions by L. Schwartz [4].

Transforming the integral expression in one variable using the Dirac delta function $\delta$ into one using the Heaviside step function $\sigma$,

$$\int_{-\infty}^{\infty} R(x)\delta(x)dx = R(0) = R(\infty) - \int_{-\infty}^{\infty} \frac{dR(x)}{dx}\sigma(x)dx = R(-\infty) + \int_{-\infty}^{\infty} \frac{dR(x)}{dx}\sigma(-x)dx \tag{1.2}$$

is essentially described in Dirac's quantum mechanics text [1]. It is derived from a well-known relation between the Dirac delta function and the derivative of the Heaviside step function. As these two expressions, the left- and right-hand sides of (1.2), are mathematically equivalent, the step-function expression would be expected to find broad applications in many branches of mathematical physics. However, that has not been noticed or developed so far compared with the application of the delta-function expression. Here we demonstrate a unified formula that extends this step-function expression for single-variable functions to multiple-variable functions. It can be interpreted as the decomposition of any nonlinear multivariate function with respect to the Heaviside step function.

## 2. Decomposition of nonlinear multivariate functions using the Heaviside step function

***Definition.*** Let $R(X_1, X_2, \cdots, X_N)$ be a continuous real function defined for $0 \leq X_i < \infty$ and satisfies:

- Whose derivatives $\frac{\partial^\alpha R}{\partial X_1^{\alpha_1} \cdots \partial X_N^{\alpha_N}}$ exist and continuous where $\alpha = \alpha_1 + \cdots + \alpha_N$ and $\alpha_i \geq 0$ is a natural number.
- $\frac{\partial^\alpha R}{\partial X_1^{\alpha_1} \cdots \partial X_N^{\alpha_N}}$ can be integrated with respect to a given $X_i$ while the other variables are held fixed.
- Sequences of functions $\left(\frac{1}{h}\left\{\frac{\partial^{\alpha-1} R}{\cdots \partial X_i^{\alpha_i - 1} \cdots}(\cdots, X_i + h, \cdots) - \frac{\partial^{\alpha-1} R}{\cdots \partial X_i^{\alpha_i - 1} \cdots}(\cdots, X_i, \cdots)\right\}\right)_{h \to 0}$ uniformly converge to $\frac{\partial^\alpha R}{\partial X_1^{\alpha_1} \cdots \partial X_N^{\alpha_N}}$ for $X_j$ and $j \neq i$.

***Theorem.*** $R(X_1, X_2, \cdots, X_N)$ can be decomposed into:

$$\begin{aligned}
R(X_1, X_2, \cdots, X_N) = R(0, 0, \cdots, 0) &+ \int_0^\infty \frac{\partial R(\mu_1, 0, \cdots, 0)}{\partial \mu_1}\sigma(X_1 - \mu_1)d\mu_1 \\
&+ \int_0^\infty \frac{\partial R(0, \mu_2, 0, \cdots, 0)}{\partial \mu_2}\sigma(X_2 - \mu_2)d\mu_2 + \cdots \\
&+ \int_0^\infty \frac{\partial R(0, \cdots, 0, \mu_i, 0, \cdots, 0)}{\partial \mu_i}\sigma(X_i - \mu_i)d\mu_i \\
&+ \cdots + \int_0^\infty \frac{\partial R(0, \cdots, 0, \mu_N)}{\partial \mu_N}\sigma(X_N - \mu_N)d\mu_N \\
&+ \int_0^\infty \int_0^\infty \frac{\partial^2 R(\mu_1, \mu_2, 0, \cdots, 0)}{\partial \mu_1 \partial \mu_2}\sigma(X_1 - \mu_1)\sigma(X_2 - \mu_2)d\mu_1\, d\mu_2 + \cdots \\
&+ \int_0^\infty \int_0^\infty \frac{\partial^2 R(0, \cdots, 0, \mu_j, 0, \cdots, 0, \mu_k, 0, \cdots, 0)}{\partial \mu_j \partial \mu_k}\sigma(X_j - \mu_j)\sigma(X_k - \mu_k)d\mu_j\, d\mu_k + \cdots
\end{aligned}$$

$$+\int_0^\infty\int_0^\infty \frac{\partial^2 R(0,\cdots,0,\mu_{N-1},\mu_N)}{\partial\mu_{N-1}\partial\mu_N}\sigma(X_{N-1}-\mu_{N-1})\sigma(X_N-\mu_N)d\mu_{N-1}\,d\mu_N$$
$$+\int_0^\infty\int_0^\infty\int_0^\infty \frac{\partial^3 R(\mu_1,\mu_2,\mu_3,0,\cdots,0)}{\partial\mu_1\partial\mu_2\partial\mu_3}\sigma(X_1-\mu_1)\sigma(X_2-\mu_2)\sigma(X_3-\mu_3)d\mu_1\,d\mu_2\,d\mu_3+\cdots$$
$$+\int_0^\infty\int_0^\infty\int_0^\infty \frac{\partial^3 R(0,\cdots,0,\mu_l,0,\cdots,0,\mu_m,0,\cdots,0,\mu_n,0,\cdots,0)}{\partial\mu_l\partial\mu_m\partial\mu_n}\sigma(X_l-\mu_l)\sigma(X_m-\mu_m)\sigma(X_n-\mu_n)d\mu_l\,d\mu_m\,d\mu_n+\cdots$$
$$+\int_0^\infty\int_0^\infty\int_0^\infty\int_0^\infty \frac{\partial^4 R(0,\cdots,0,\mu_o,0,\cdots,0,\mu_p,0,\cdots,0,\mu_q,0,\cdots,0,\mu_r,0,\cdots,0)}{\partial\mu_o\partial\mu_p\partial\mu_q\partial\mu_r}\sigma(X_o-\mu_o)\sigma(X_p-\mu_p)\sigma(X_q-\mu_q)\sigma(X_r-\mu_r)\,d\mu_o\,d\mu_p\,d\mu_q\,d\mu_r+\cdots$$
$$+\int_0^\infty\cdots\int_0^\infty \frac{\partial^N R(\mu_1,\cdots,\mu_N)}{\partial\mu_1\cdots\partial\mu_N}\sigma(X_1-\mu_1)\cdots\sigma(X_N-\mu_N)d\mu_1\cdots d\mu_N\,, \tag{2.1}$$

where

$$\sigma(X_i-\mu_i)=\begin{cases}1 & (X_i>\mu_i)\\ \frac{1}{2} & (X_i=\mu_i)\\ 0 & (X_i<\mu_i)\end{cases} \tag{2.2}$$

defines a set of Heaviside step functions for each $\mu_i \geq 0$.

**Proof.** This proof is demonstrated using mathematical induction.

Using the definition of the Dirac delta function, then for any real function $R(X_1,X_2,\cdots,X_N)$

$$R(X_1,X_2,\cdots,X_N)=\int_{-\infty}^{\infty}\cdots\int_{-\infty}^{\infty}R(\mu_1,\cdots,\mu_N)\delta(\mu_1-X_1)\cdots\delta(\mu_N-X_N)d\mu_1\cdots d\mu_N\,. \tag{2.3}$$

Therefore,

$$R(X_1,X_2,\cdots,X_N)=\int_0^\infty\cdots\int_0^\infty R(\mu_1,\cdots,\mu_N)\delta(\mu_1-X_1)\cdots\delta(\mu_N-X_N)d\mu_1\cdots d\mu_N \tag{2.4}$$

with $X_i \geq 0$.

For single-variable functions ($N=1$),

$$\int_0^\infty R(\mu_1)\delta(\mu_1-X_1)d\mu_1=R(0)+\int_0^\infty \frac{d}{d\mu_1}R(\mu_1)\sigma(X_1-\mu_1)d\mu_1 \tag{2.5}$$

holds by Lemma 2.1 given below. Note that, for single-variable functions, Dirac described the essentially equivalent expression (1.2) [1]. To initiate the mathematical induction procedure, suppose that

$$\int_0^\infty\cdots\int_0^\infty R(\mu_1,\cdots,\mu_N)\delta(\mu_1-X_1)\cdots\delta(\mu_N-X_N)d\mu_1\cdots d\mu_N$$
$$=R(0,0,\cdots,0)+\int_0^\infty \frac{\partial R(\mu_1,0,\cdots,0)}{\partial\mu_1}\sigma(X_1-\mu_1)d\mu_1+\cdots$$
$$+\int_0^\infty \frac{\partial R(0,\cdots,0,\mu_i,0,\cdots,0)}{\partial\mu_i}\sigma(X_i-\mu_i)d\mu_i+\cdots$$
$$+\int_0^\infty\int_0^\infty \frac{\partial^2 R(0,\cdots,0,\mu_j,0,\cdots,0,\mu_k,0,\cdots,0)}{\partial\mu_j\partial\mu_k}\sigma(X_j-\mu_j)\sigma(X_k-\mu_k)d\mu_j\,d\mu_k+\cdots$$
$$+\int_0^\infty\cdots\int_0^\infty \frac{\partial^N R(\mu_1,\cdots,\mu_N)}{\partial\mu_1\cdots\partial\mu_N}\sigma(X_1-\mu_1)\cdots\sigma(X_N-\mu_N)d\mu_1\cdots d\mu_N \tag{2.6}$$

holds for some *N*. Using (2.6), the following,

$$\int_0^\infty\cdots\int_0^\infty R(\mu_1,\cdots,\mu_N,\mu_{N+1})\delta(\mu_1-X_1)\cdots\delta(\mu_N-X_N)d\mu_1\cdots d\mu_N$$
$$=R(X_1,\cdots,X_N,\mu_{N+1})$$
$$=R(0,0,\cdots,0,\mu_{N+1})+\cdots$$
$$+\int_0^\infty \frac{\partial R(0,\cdots,0,\mu_i,0,\cdots,0,\mu_{N+1})}{\partial\mu_i}\sigma(X_i-\mu_i)d\mu_i+\cdots$$
$$+\int_0^\infty\cdots\int_0^\infty \frac{\partial^N R(\mu_1,\cdots,\mu_N,\mu_{N+1})}{\partial\mu_1\cdots\partial\mu_N}\sigma(X_1-\mu_1)\cdots\sigma(X_N-\mu_N)d\mu_1\cdots d\mu_N\,, \tag{2.7}$$

holds because $R(X_1,\cdots,X_N,\mu_{N+1})$ can be regarded as one of the $R(X_1,\cdots,X_N)$ appending a parameter $\mu_{N+1}$. Multiplying both sides of (2.7) by $\delta(\mu_{N+1}-X_{N+1})$ and then integrating each term with respect to $\mu_{N+1}$, one obtains on the left-hand side,

$$\int_0^\infty\cdots\int_0^\infty R(\mu_1,\cdots,\mu_N,\mu_{N+1})\delta(\mu_1-X_1)\cdots\delta(\mu_N-X_N)\delta(\mu_{N+1}-X_{N+1})d\mu_1\cdots d\mu_N d\mu_{N+1}\,, \tag{2.8}$$

and on the right-hand side,

$$\int_0^\infty\{R(0,0,\cdots,0,\mu_{N+1})\}\delta(\mu_{N+1}-X_{N+1})d\mu_{N+1}$$
$$+\cdots+\int_0^\infty\left\{\int_0^\infty \frac{\partial R(0,\cdots,0,\mu_i,0,\cdots,0,\mu_{N+1})}{\partial\mu_i}\sigma(X_i-\mu_i)d\mu_i\right\}\delta(\mu_{N+1}-X_{N+1})d\mu_{N+1}$$
$$+\cdots+\int_0^\infty\left\{\int_0^\infty\cdots\int_0^\infty \frac{\partial^N R(\mu_1,\cdots,\mu_N,\mu_{N+1})}{\partial\mu_1\cdots\partial\mu_N}\sigma(X_1-\mu_1)\cdots\sigma(X_N-\mu_N)d\mu_1\cdots d\mu_N\right\}\delta(\mu_{N+1}-X_{N+1})d\mu_{N+1}\,. \tag{2.9}$$

The order of the integrations can be changed because each integrand can be integrated with respect to its corresponding $\mu_i$ while holding other variables fixed. Therefore, (2.9) can be transformed into

$$\left\{\int_0^\infty R(0,0,\cdots,0,\mu_{N+1})\delta(\mu_{N+1}-X_{N+1})d\mu_{N+1}\right\}$$
$$+\cdots+\int_0^\infty\left\{\int_0^\infty \frac{\partial R(0,\cdots,0,\mu_i,0,\cdots,0,\mu_{N+1})}{\partial\mu_i}\delta(\mu_{N+1}-X_{N+1})d\mu_{N+1}\right\}\sigma(X_i-\mu_i)d\mu_i$$

$$+\cdots+\int_0^\infty\cdots\int_0^\infty\left\{\int_0^\infty\frac{\partial^N R(\mu_1,\cdots,\mu_N,\mu_{N+1})}{\partial\mu_1\cdots\partial\mu_N}\delta(\mu_{N+1}-X_{N+1})d\mu_{N+1}\right\}\sigma(X_1-\mu_1)\cdots\sigma(X_N-\mu_N)d\mu_1\cdots d\mu_N\,. \tag{2.10}$$

The terms enclosed in braces in (2.10) can be transformed using Lemma 2.1 as follows,

$$\left\{R(0,0,\cdots,0,0)+\int_0^\infty\frac{\partial R(0,\cdots,0,\mu_{N+1})}{\partial\mu_{N+1}}\sigma(X_{N+1}-\mu_{N+1})d\mu_{N+1}\right\}$$
$$+\cdots+\int_0^\infty\left\{\frac{\partial R(0,\cdots,0,\mu_i,0,\cdots,0,0)}{\partial\mu_i}+\int_0^\infty\frac{\partial^2 R(0,\cdots,0,\mu_i,0,\cdots,0,\mu_{N+1})}{\partial\mu_{N+1}\partial\mu_i}\sigma(X_{N+1}-\mu_{N+1})d\mu_{N+1}\right\}\sigma(X_i-\mu_i)d\mu_i$$
$$+\cdots+\int_0^\infty\cdots\int_0^\infty\left\{\frac{\partial^N R(\mu_1,\cdots,\mu_N,0)}{\partial\mu_1\cdots\partial\mu_N}+\int_0^\infty\frac{\partial^{N+1} R(\mu_1,\cdots,\mu_N,\mu_{N+1})}{\partial\mu_{N+1}\partial\mu_1\cdots\partial\mu_N}\sigma(X_{N+1}-\mu_{N+1})d\mu_{N+1}\right\}\sigma(X_1-\mu_1)\cdots\sigma(X_N-\mu_N)d\mu_1\cdots d\mu_N\,, \tag{2.11}$$

where Lemma 2.4 was also used.

Finally (2.11) becomes

$$R(0,0,\cdots,0,0)$$
$$+\int_0^\infty\frac{\partial R(0,\cdots,0,\mu_{N+1})}{\partial\mu_{N+1}}\sigma(X_{N+1}-\mu_{N+1})d\mu_{N+1}+\cdots$$
$$+\int_0^\infty\frac{\partial R(0,\cdots,0,\mu_i,0,\cdots,0,0)}{\partial\mu_i}\sigma(X_i-\mu_i)d\mu_i$$
$$+\cdots+\int_0^\infty\int_0^\infty\frac{\partial^2 R(0,\cdots,0,\mu_i,0,\cdots,0,\mu_{N+1})}{\partial\mu_i\partial\mu_{N+1}}\sigma(X_i-\mu_i)\sigma(X_{N+1}-\mu_{N+1})d\mu_i d\mu_{N+1}$$
$$+\cdots+\int_0^\infty\cdots\int_0^\infty\frac{\partial^N R(\mu_1,\cdots,\mu_N,0)}{\partial\mu_1\cdots\partial\mu_N}\sigma(X_1-\mu_1)\cdots\sigma(X_N-\mu_N)d\mu_1\cdots d\mu_N+\cdots$$
$$+\int_0^\infty\cdots\int_0^\infty\frac{\partial^{N+1} R(\mu_1,\cdots,\mu_N,\mu_{N+1})}{\partial\mu_{N+1}\partial\mu_1\cdots\partial\mu_N}\sigma(X_1-\mu_1)\cdots\sigma(X_N-\mu_N)\sigma(X_{N+1}-\mu_{N+1})d\mu_1\cdots d\mu_N d\mu_{N+1}\,. \tag{2.12}$$

Thus, assuming expression (2.6) for $N$ leads to the same expression for $N+1$. The formula for $N=1$ also holds as described above. Therefore, the theorem holds for any natural number $N$. □

***Lemma 2.1.***

$$\int_0^\infty F(\mu_1,\cdots,\mu_i,\cdots,\mu_N)\delta(\mu_i-X_i)d\mu_i$$
$$=F(\mu_1,\cdots,\mu_{i-1},0,\mu_{i+1},\cdots,\mu_N)+\int_0^\infty\frac{\partial F(\mu_1,\cdots,\mu_i,\cdots,\mu_N)}{\partial\mu_i}\sigma(X_i-\mu_i)d\mu_i\,. \tag{2.13}$$

**Proof.** From Lemma 2.2,

$$\int_0^\infty F(\mu_1,\cdots,\mu_i,\cdots,\mu_N)\delta(\mu_i-X_i)d\mu_i$$
$$=\int_0^\infty F(\mu_1,\cdots,\mu_i,\cdots,\mu_N)\frac{d}{d\mu_i}\sigma(\mu_i-X_i)d\mu_i$$
$$=-\int_0^\infty F(\mu_1,\cdots,\mu_i,\cdots,\mu_N)\frac{d}{d\mu_i}\sigma(X_i-\mu_i)d\mu_i\,. \tag{2.14}$$

If $X_i>0$,

$$-\int_0^\infty F(\mu_1,\cdots,\mu_i,\cdots,\mu_N)\frac{d}{d\mu_i}\sigma(X_i-\mu_i)d\mu_i$$
$$=-[F(\mu_1,\cdots,\mu_i,\cdots,\mu_N)\sigma(X_i-\mu_i)]_0^\infty+\int_0^\infty\frac{\partial}{\partial\mu_i}F(\mu_1,\cdots,\mu_i,\cdots,\mu_N)\sigma(X_i-\mu_i)d\mu_i$$
$$=-F(\mu_1,\cdots,\infty,\cdots,\mu_N)\sigma(X_i-\infty)+F(\mu_1,\cdots,0,\cdots,\mu_N)\sigma(X_i-0)+\int_0^\infty\frac{\partial}{\partial\mu_i}F(\mu_1,\cdots,\mu_i,\cdots,\mu_N)\sigma(X_i-\mu_i)d\mu_i\,. \tag{2.15}$$

With $\sigma(X_i-\infty)=0$ and $\sigma(X_i-0)=1$, we obtain

$$\int_0^\infty F(\mu_1,\cdots,\mu_i,\cdots,\mu_N)\delta(\mu_i-X_i)d\mu_i$$
$$=F(\mu_1,\cdots,0,\cdots,\mu_N)+\int_0^\infty\frac{\partial F(\mu_1,\cdots,\mu_i,\cdots,\mu_N)}{\partial\mu_i}\sigma(X_i-\mu_i)d\mu_i\,. \tag{2.16}$$

If $X_i=0$, the left-hand side of (2.13) is

$$\int_0^\infty F(\mu_1,\cdots,\mu_i,\cdots,\mu_N)\delta(\mu_i-X_i)d\mu_i$$
$$=\int_0^\infty F(\mu_1,\cdots,\mu_i,\cdots,\mu_N)\delta(\mu_i)d\mu_i$$
$$=F(\mu_1,\cdots,\mu_{i-1},0,\mu_{i+1},\cdots,\mu_N) \tag{2.17}$$

whereas the right-hand side of (2.13) is

$$F(\mu_1,\cdots,\mu_{i-1},0,\mu_{i+1},\cdots,\mu_N)$$
$$+\int_0^\infty\frac{\partial F(\mu_1,\cdots,\mu_i,\cdots,\mu_N)}{\partial\mu_i}\sigma(X_i-\mu_i)d\mu_i$$

$$= F(\mu_1, \cdots, \mu_{i-1}, 0, \mu_{i+1}, \cdots, \mu_N) + \int_0^0 \frac{\partial F(\mu_1, \cdots, \mu_i, \cdots, \mu_N)}{\partial \mu_i} \sigma(X_i - \mu_i) d\mu_i + \int_{0<}^{\infty} \frac{\partial F(\mu_1, \cdots, \mu_i, \cdots, \mu_N)}{\partial \mu_i} \sigma(X_i - \mu_i) d\mu_i = F(\mu_1, \cdots, \mu_{i-1}, 0, \mu_{i+1}, \cdots, \mu_N) . \tag{2.18}$$

□

***Lemma 2.2.***

$$\frac{d\sigma(\mu - X)}{d\mu} = -\frac{d\sigma(X - \mu)}{d\mu} . \tag{2.19}$$

**Proof.** From Lemma 2.3, $\sigma(\mu - X) = 1 - \sigma(X - \mu)$ . Therefore, by differentiating both sides, $\frac{d\sigma(\mu-X)}{d\mu} = -\frac{d\sigma(X-\mu)}{d\mu}$. □

***Lemma 2.3.***

$$\sigma(-X) = 1 - \sigma(X) . \tag{2.20}$$

**Proof.** If $X > 0$, $\sigma(-X) = 0 = 1 - \sigma(X)$ . If $X < 0$, $\sigma(-X) = 1 = 1 - \sigma(X)$ . If $X = 0$, $\sigma(-X) = \frac{1}{2} = 1 - \sigma(X)$ . □

***Lemma 2.4.***

$$\left. \frac{\partial^\alpha R(\cdots, \mu_i, \cdots, \mu_j, \cdots)}{\cdots \partial \mu_i^{\alpha_i} \cdots} \right|_{\substack{\mu_j = 0 \\ (i \neq j)}} = \frac{\partial^\alpha R(\cdots, \mu_i, \cdots, \mu_{j-1}, 0, \mu_{j+1}, \cdots)}{\cdots \partial \mu_i^{\alpha_i} \cdots} , \tag{2.21}$$

where $\alpha = \alpha_1 + \cdots + \alpha_N$ and $\alpha_i \geq 0$ is a natural number.

**Proof.** Since derivatives of $R$ is continuous,

$$\left. \frac{\partial^\alpha R(\cdots, \mu_i, \cdots, \mu_j, \cdots)}{\cdots \partial \mu_i^{\alpha_i} \cdots} \right|_{\substack{\mu_j = 0 \\ (i \neq j)}} = \lim_{\substack{\mu_j \to 0 \\ (i \neq j)}} \frac{\partial R^\alpha(\cdots, \mu_i, \cdots, \mu_j, \cdots)}{\cdots \partial \mu_i^{\alpha_i} \cdots} = \lim_{\substack{\mu_j \to 0 \\ (i \neq j)}} \lim_{h \to 0} \frac{1}{h} \left\{ \frac{\partial^{\alpha-1} R}{\cdots \partial \mu_i^{\alpha_i - 1} \cdots} (\cdots, \mu_i + h, \cdots) - \frac{\partial^{\alpha-1} R}{\cdots \partial \mu_i^{\alpha_i - 1} \cdots} (\cdots, \mu_i, \cdots) \right\} = \lim_{\substack{\mu_j \to 0 \\ (i \neq j)}} (r_{h \to 0}), \tag{2.22}$$

where $(r_{h \to 0})$ is a sequence of functions. The sequence of function $(r_{h \to 0})$ converges uniformly for $\mu_j$ and $j \neq i$ to $\frac{\partial^\alpha R(\cdots, \mu_i, \cdots, \mu_j, \cdots)}{\cdots \partial \mu_i^{\alpha_i} \cdots}$. Similarly $\lim_{\substack{\mu_j \to 0 \\ (i \neq j)}} \frac{1}{h} \left\{ \frac{\partial^{\alpha-1} R}{\cdots \partial \mu_i^{\alpha_i - 1} \cdots} (\cdots, \mu_i + h, \cdots, \mu_j, \cdots) - \frac{\partial^{\alpha-1} R}{\cdots \partial \mu_i^{\alpha_i - 1} \cdots} (\cdots, \mu_i, \cdots, \mu_j, \cdots) \right\}$ converges pointwise for $h$ to $\frac{1}{h} \left\{ \frac{\partial^{\alpha-1} R}{\cdots \partial \mu_i^{\alpha_i - 1} \cdots} (\cdots, \mu_i + h, \cdots, \mu_j = 0, \cdots) - \frac{\partial^{\alpha-1} R}{\cdots \partial \mu_i^{\alpha_i - 1} \cdots} (\cdots, \mu_i, \cdots, \mu_j = 0, \cdots) \right\}$ since it is continuous. Therefore, the order of the limits can be interchanged. Finally (2.22) be

$$\lim_{\substack{\mu_j \to 0 \\ (i \neq j)}} (r_{h \to 0}) = \lim_{h \to 0} \lim_{\substack{\mu_j \to 0 \\ (i \neq j)}} \frac{1}{h} \left\{ \frac{\partial^{\alpha-1} R}{\cdots \partial \mu_i^{\alpha_i - 1} \cdots} (\cdots, \mu_i + h, \cdots, \mu_j, \cdots) - \frac{\partial^{\alpha-1} R}{\cdots \partial \mu_i^{\alpha_i - 1} \cdots} (\cdots, \mu_i, \cdots, \mu_j, \cdots) \right\} = \lim_{h \to 0} \frac{1}{h} \left\{ \frac{\partial^{\alpha-1} R}{\cdots \partial \mu_i^{\alpha_i - 1} \cdots} (\cdots, \mu_i + h, \cdots, \mu_j = 0, \cdots) - \frac{\partial^{\alpha-1} R}{\cdots \partial \mu_i^{\alpha_i - 1} \cdots} (\cdots, \mu_i, \cdots, \mu_j = 0, \cdots) \right\} = \frac{\partial^\alpha R(\cdots, \mu_i, \cdots, \mu_{j-1}, 0, \mu_{j+1}, \cdots)}{\cdots \partial \mu_i^{\alpha_i} \cdots}. \tag{2.23}$$

□

## 3. Concluding remarks

We have demonstrated the decomposition of a nonlinear multivariate function as a sum of integrals of which each integrand is composed of a derivative and a direct product of Heaviside step functions. The expression offers a mathematically unified and systematic expansion equivalent to that given in terms of Dirac delta functions. It can further be approximated using, e.g., sigmoid functions with suitable parameters, yielding a convenient form for broad applications to other fields exploiting both analytical and numerical methods.

**Acknowledgements:**

I thank K. Soda for discussion.